\documentclass[12pt]{article}
\usepackage{amsmath}
\usepackage{amssymb}
\usepackage{tikz}
\newtheorem{theorem}{Theorem}
\newtheorem{corollary}{Corollary}

\newtheorem{proposition}{Proposition}
\newtheorem{question}{Question}
\begin{document}
\title{On $2$-colored graphs and partitions of boxes}
\author{Ron Holzman\thanks{Department of Mathematics, Technion-Israel Institute of Technology, 32000 Haifa, Israel. E-mail: holzman@technion.ac.il. Research supported in part by ISF grant no. 1162/15. Helpful discussions with Joseph Briggs are gratefully acknowledged.}}
\maketitle

\begin{abstract}
We prove that if the edges of a graph $G$ can be colored blue or red in such a way that every vertex belongs to a monochromatic $k$-clique of each color, then $G$ has at least $4(k-1)$ vertices. This confirms a conjecture of Bucic et al. \cite{BLLW}, and thereby solves the $2$-dimensional case of their problem about partitions of discrete boxes with the $k$-piercing property. We also characterize the case of equality in our result.
\end{abstract}

\section{Introduction}

In this paper, a \emph{$2$-colored graph} will be a simple graph with edges colored blue or red. Bucic et al. \cite{BLLW} asked the following: Given an integer $k \ge 2$, what is the smallest possible number of vertices in a $2$-colored graph having the property that every vertex belongs to a monochromatic $k$-clique of each color?

They gave the following construction, showing that $4(k-1)$ vertices suffice. First, for $k=2$, take a $4$-cycle with edges colored alternatingly. Now, for general $k$, blow up this $4$-cycle, replacing each vertex by a monochromatic $(k-1)$-clique, with colors alternating along the $4$-cycle (all edges between two adjacent $(k-1)$-cliques get the same color as the edge in the underlying $4$-cycle). It is easy to verify that this $2$-colored graph has the required property.

Bucic et al. \cite{BLLW} conjectured that this construction is optimal, and proved a lower bound of the form $(4 - o_k(1))k$ on the number of vertices in any $2$-colored graph with the required property. Here we prove exact optimality.

\begin{theorem} \label{main}
Let $k \ge 2$ be an integer, and let $G=(V,E)$ be a $2$-colored graph so that every vertex in $V$ belongs to a monochromatic $k$-clique of each color. Then $|V| \ge 4(k-1)$.
\end{theorem}

Our proof, given in Section~2, combines counting arguments with a linear algebraic trick similar to one used by Tverberg \cite{T}. In Section~3 we characterize the case of equality in Theorem~\ref{main}. Perhaps surprisingly, for $k \ge 3$ it turns out that the above example on $4(k-1)$ vertices is not the only extremal one. In Section~4 we discuss some generalizations and reformulations of the above question. These involve, in particular, partitions of a box into sub-boxes and decompositions of a bipartite graph into complete bipartite subgraphs.

\section{Proof of Theorem~\ref{main}}

Let $G=(V,E)$ be a $2$-colored graph so that every vertex in $V$ belongs to a monochromatic $k$-clique of each color. Instead of working directly with the graph, we store the information we have in the following form: We have a vertex-set $V$ and two families $\mathcal{B}=\{B_1,\ldots,B_b\}$, $\mathcal{R}=\{R_1,\ldots,R_r\}$ of subsets of $V$ satisfying:

\begin{equation} \label{one}
|B_i| \ge k \textrm{ and } |R_j| \ge k \quad \forall i \in [b], j \in [r],
\end{equation}

\begin{equation} \label{two}
|B_i \cap R_j| \le 1 \quad \forall i \in [b], j \in [r],
\end{equation}

\begin{equation} \label{three}
\bigcup_{i=1}^b B_i = \bigcup_{j=1}^r R_j = V.
\end{equation}

\noindent Indeed, given $G$ we can construct $\mathcal{B}$ and $\mathcal{R}$ as the families of vertex-sets of blue (resp. red) cliques witnessing that every vertex belongs to a monochromatic clique of each color.

Next, we claim that we can keep the same vertex-set $V$ and possibly make some adjustments to the families $\mathcal{B}$ and $\mathcal{R}$, so that the following will be satisfied in addition to (\ref{one})--(\ref{three}):

\begin{equation} \label{four}
B_i \setminus \bigcup_{i' \ne i} B_{i'} \ne \emptyset \textrm{ and } R_j \setminus \bigcup_{j' \ne j} R_{j'} \ne \emptyset \quad \forall i \in [b], j \in [r],
\end{equation}

\begin{equation} \label{five}
|B_i \cap R_j| = 1 \quad \forall i \in [b], j \in [r].
\end{equation}

\noindent Indeed, if $B_i \subseteq \bigcup_{i' \ne i} B_{i'}$, say, then we can discard $B_i$ from $\mathcal{B}$ while retaining properties (\ref{one})--(\ref{three}). Iterating this operation we end up with families satisfying (\ref{one})--(\ref{four}). At this point, if $B_i \cap R_j = \emptyset$ then we can choose a vertex $v \in B_i \setminus \bigcup_{i' \ne i} B_{i'}$ and replace $R_j$ by $R_j \cup \{v\}$ while retaining properties (\ref{one})--(\ref{three}). It may happen that this change causes a violation of (\ref{four}), namely when we had $R_{j^*} \setminus \bigcup_{j' \ne j^*} R_{j'} = \{v\}$ for some $j^* \ne j$ before the change; in this case, after adding $v$ to $R_j$ we discard $R_{j^*}$. Iterating this operation we end up with families satisfying (\ref{one})--(\ref{five}).

Thus, we may assume that the set $V$ and the two families $\mathcal{B}=\{B_1,\ldots,B_b\}$, $\mathcal{R}=\{R_1,\ldots,R_r\}$ of subsets of $V$ satisfy (\ref{one})--(\ref{five}). For every $v \in V$ we write

\vspace{-10pt}

\begin{equation*}
I_v = \{ i \in [b]: v \in B_i \}, \quad J_v = \{ j \in [r]: v \in R_j \}.
\end{equation*}

\noindent Note that the properties of $V$, $\mathcal{B}$ and $\mathcal{R}$ may be expressed in terms of the subsets $I_v$ of $[b]$ and $J_v$ of $[r]$, for $v \in V$, as follows: (\ref{one}) says that every $i \in [b]$ is covered by the subsets $I_v$ at least $k$ times, and similarly for $[r]$ and the subsets $J_v$; (\ref{three}) says that $I_v, J_v \ne \emptyset$; (\ref{four}) says that for every $i \in [b]$ there is $v \in V$ such that $I_v = \{i\}$, and similarly for $[r]$ and $J_v$; (\ref{five}) says that the product sets $I_v \times J_v$ partition $[b] \times [r]$.

\begin{proposition} \label{prop}
If $V$, $\mathcal{B}$ and $\mathcal{R}$ satisfy (\ref{four}) and (\ref{five}) then $|V| \ge b + r - 1$.
\end{proposition}

\noindent \textbf{Proof} We introduce for each $i \in [b]$ a variable $x_i$, and for each $j \in [r]$ a variable $y_j$ (these variables take real values). By (\ref{five}) we have the identity

\begin{equation} \label{six}
\sum_{v \in v} (\sum_{i \in I_v} x_i)(\sum_{j \in J_v} y_j) = (\sum_{i=1}^b x_i)(\sum_{j=1}^r y_j).
\end{equation}

\noindent Now we consider the following system of homogeneous linear equations:

\begin{equation} \label{seven}
\sum_{i \in I_v} x_i - \sum_{j \in J_v} y_j = 0, \quad v \in V,
\end{equation}

\begin{equation} \label{eight}
\sum_{i=1}^b x_i = 0.
\end{equation}

\noindent It suffices to show that the system has only the trivial solution, because this implies that the number of equations $|V|+1$ is at least as large as the number of variables $b+r$. Let $(x_i)_{i \in [b]}, (y_j)_{j \in [r]}$ satisfy (\ref{seven}) and (\ref{eight}). By (\ref{seven}) we know that for each $v \in V$ there is a real number $\alpha_v$ so that $\sum_{i \in I_v} x_i = \sum_{j \in J_v} y_j = \alpha_v$. The identity (\ref{six}) implies, using (\ref{eight}), that $\sum_{v \in V} \alpha_v^2 = 0$ and hence $\alpha_v = 0$ for all $v \in V$. Now, given $i \in [b]$ we can find by (\ref{four}) some $v \in V$ such that $x_i = \sum_{i \in I_v} x_i = \alpha_v = 0$, and a similar argument shows that $y_j=0$ for every $j \in [r]$, as required. $\square$

\bigskip

Returning to the proof of Theorem~\ref{main}, we may henceforth assume that $b+r \le 4(k-1)$, otherwise $|V| \ge 4(k-1)$ follows from Proposition~\ref{prop}. We also know that $b \ge k$, because the sets $I_v$, $v \in R_1$, are $k$ or more disjoint nonempty subsets of $[b]$; similarly $r \ge k$. Thus, the relevant domain for $b+r$ in the rest of the proof is

\begin{equation} \label{nine}
2k \le b+r \le 4(k-1).
\end{equation}

Using (\ref{one}) we have

\begin{equation} \label{ten}
\sum_{v \in V} |I_v| + |J_v| \ge k(b+r),
\end{equation}

\noindent and using (\ref{five}) we have

\begin{equation} \label{eleven}
\sum_{v \in V} |I_v||J_v| = br.
\end{equation}

\noindent Since $|I_v|$ and $|J_v|$ are nonzero by (\ref{three}), their product is smallest (given their sum) when one of them is $1$. Hence

\begin{equation} \label{twelve}
|I_v||J_v| \ge |I_v| + |J_v| - 1 \quad \forall v \in V.
\end{equation}

\noindent Using (\ref{ten})--(\ref{twelve}) we can write

\vspace{-15pt}

\begin{eqnarray}
|V| & = & \sum_{v \in V} |I_v| + |J_v| - (|I_v| + |J_v| - 1) \nonumber\\
 & \ge & k(b+r) - \sum_{v \in V} (|I_v| + |J_v| - 1) \nonumber \\
 & \ge & k(b+r) - \sum_{v \in V} |I_v||J_v| \label{thirteen} \\
 & = & k(b+r) - br \nonumber \\
 & \ge & k(b+r) - \frac{(b+r)^2}{4}. \nonumber
\end{eqnarray}

\noindent The latter is a decreasing function of $b+r$ in the domain (\ref{nine}), and is therefore bounded from below by its value at $b+r = 4(k-1)$, which is $4(k-1)$. This proves that $|V| \ge 4(k-1)$, as required. $\square$

\section{Characterization of extremal graphs}

If $G=(V,E)$ is a $2$-colored graph having the property that every vertex in $V$ belongs to a monochromatic $k$-clique of each color, then adding any edges to $G$ (between existing vertices) and coloring them arbitrarily results in a graph with the same property. Therefore we can restrict attention to those graphs having this property which are \emph{edge-critical}, in the sense that removing any edge entails the loss of this property.

Here is a construction of an edge-critical $2$-colored graph on $4(k-1)$ vertices, so that every vertex belongs to a monochromatic $k$-clique of each color, which generalizes the one from~\cite{BLLW} described in the introduction. Let $k \ge 2$ be an integer, let $X$ and $Y$ be two disjoint sets of $2(k-1)$ vertices each, and let $B(X,Y)$ and $R(X,Y)$ be two complementary $(k-1)$-regular bipartite graphs on the bipartition $(X,Y)$. Our graph $G=G(X,Y,B,R)$ has $X \cup Y$ as its vertex-set. It has the complete bipartite graph on $(X,Y)$ as a subgraph, with edges in $B(X,Y)$ colored blue and edges in $R(X,Y)$ colored red. We refer to $B(X,Y)$ and $R(X,Y)$ as the blue and red graphs, respectively. In addition, any two vertices in $X$ which have a common neighbor in the blue graph are joined by a blue edge in $G$, and any two vertices in $Y$ which have a common neighbor in the red graph are joined by a red edge in $G$. It is easy to verify that this $2$-colored graph has the required property and is edge-critical.

For $k=2$ we have $|X|=|Y|=2$ and the blue and red graphs must be two complementary perfect matchings, resulting in the $2$-colored $4$-cycle described in the introduction. But for higher values of $k$, we have more freedom in choosing $B(X,Y)$ and $R(X,Y)$. For example, consider $k=3$, so $|X|=|Y|=4$. We may choose $B(X,Y)$ and $R(X,Y)$ so that each of them is the disjoint union of two $4$-cycles, resulting in the blown-up $4$-cycle graph from the introduction. But we can also choose $B(X,Y)$ and $R(X,Y)$ to be $8$-cycles, resulting in a new example, not isomorphic to the previous one.

Note that the construction described in the introduction corresponds to the following choice of $B(X,Y)$ and $R(X,Y)$: $X$ is equi-partitioned into $X_1$ and $X_2$, $Y$ is equi-partitioned into $Y_1$ and $Y_2$, $B(X,Y)$ consists of all edges between $X_1$ and $Y_1$ and between $X_2$ and $Y_2$, and $R(X,Y)$ consists of all edges between $X_1$ and $Y_2$ and between $X_2$ and $Y_1$. For this choice, the resulting graph $G(X,Y,B,R)$ induces blue cliques on $X_1$ and $X_2$ and red cliques on $Y_1$ and $Y_2$, and has a total of $2(k-1)(3k-4)$ edges. Among all graphs of the form $G(X,Y,B,R)$ for a given value of $k$, the latter uniquely minimizes the number of edges. To see this, observe that in the graph induced on $X$ (and similarly for $Y$) each vertex must have degree at least $k-2$, and the only way to have these degrees equal to $k-2$ is by using $X_1,X_2,Y_1,Y_2$ as above.

The next result shows that all edge-critical extremal examples for Theorem~\ref{main} are of the form $G=G(X,Y,B,R)$, thus characterizing the case of equality in that theorem.

\begin{theorem} \label{equality}
Let $k \ge 2$ be an integer, and let $|V|=4(k-1)$. Let $G=(V,E)$ be a $2$-colored graph so that every vertex in $V$ belongs to a monochromatic $k$-clique of each color, and $G$ is edge-critical with respect to this property. Then $G$ is isomorphic to some $G(X,Y,B,R)$, where $B(X,Y)$ and $R(X,Y)$ are complementary $(k-1)$-regular bipartite graphs on $(X,Y)$.
\end{theorem}

\noindent \textbf{Proof} Let $G=(V,E)$ satisfy the assumptions of the theorem. In the case $k=2$, it is easy to check directly that $G$ must be a $4$-cycle colored alternatingly, as claimed. We henceforth assume that $k \ge 3$.

As in the proof of Theorem~\ref{main}, we associate with $G$ two families $\mathcal{B}=\{B_1,\ldots,B_b\}$, $\mathcal{R}=\{R_1,\ldots,R_r\}$ of subsets of $V$ satisfying (\ref{one})--(\ref{three}). Clearly, the blue edges of $G$ are those pairs $\{u,v\}$ contained in some $B_i$, and the red edges are those pairs $\{u,v\}$ contained in some $R_j$ (by edge-criticality, there can be no other edges in $G$). In the main part of the proof below, we assume that $V$, $\mathcal{B}$, $\mathcal{R}$ satisfy (\ref{four}) and (\ref{five}) as well; at the end of the proof we will justify this assumption. We also use the notations $I_v$ and $J_v$ for $v \in V$ as introduced in the proof of Theorem~\ref{main}. According to that proof, the only values of $b+r$ which may result in $|V|=4(k-1)$ are $4(k-1)$ and $4(k-1)+1$ (if $b+r < 4(k-1)$ then (\ref{thirteen}) forces $|V|$ to be larger, and if $b+r > 4(k-1)+1$ then Proposition~\ref{prop} does that).

\newpage

\noindent \textbf{Case 1} $\quad b+r=4(k-1)$

\bigskip

Because $|V|=4(k-1)$, (\ref{thirteen}) must hold as an equality. This implies that (\ref{ten}) and (\ref{twelve}) hold as equalities, and $b=r=2(k-1)$. Equality in (\ref{ten}) means that every $B_i$ and every $R_j$ is of size exactly $k$. Equality in (\ref{twelve}) means that for every $v \in V$, at least one of $I_v, J_v$ is a singleton. For $j \in [r]$, the sets $I_v$, $v \in R_j$, partition $[b]$ into $k$ nonempty subsets. This implies that $|I_v| \le k-1$, and similarly $|J_v| \le k-1$, for every $v \in V$. Therefore $|I_v||J_v| \le k-1$ for every $v \in V$, but since $\sum_{v \in V} |I_v||J_v| = 4(k-1)^2$ we must have equality for every $v \in V$. This means that we can partition $V$ into two sets:

\vspace{-15pt}

\begin{equation*}
X=\{v \in V: |I_v|=k-1, |J_v|=1 \}, \quad Y=\{v \in V: |I_v|=1, |J_v|=k-1 \}.
\end{equation*}

\noindent As $\sum_{v \in V} |I_v| = kb = 2k(k-1)$, we must have $|X|=|Y|=2(k-1)$.

Now, consider a vertex $v \in X$. There is $j \in [r]$ such that $v \in R_j$. Since the sets $I_u$, $u \in R_j$, partition $[b]$ into $k$ subsets, one of which is $I_v$ of size $k-1$, all other $I_u$ must be singletons, so that $R_j \setminus \{v\} \subseteq Y$. This accounts for $k-1$ red edges from $v$ into $Y$. As this holds for every $v \in X$, and similarly every $v \in Y$ must have at least $k-1$ blue edges into $X$, the complete bipartite graph on $(X,Y)$ must appear in $G$ and be colored so that the blue graph $B(X,Y)$ and the red graph $R(X,Y)$ are both $(k-1)$-regular. The above also implies that the neighbors of every $v \in X$ in the red graph must form a red clique in $Y$, and the neighbors of every $v \in Y$ in the blue graph must form a blue clique in $X$. This shows that $G(X,Y,B,R)$ is contained in $G$, and as $G$ is edge-critical, they must coincide.

\bigskip

\noindent \textbf{Case 2} $\quad b+r=4(k-1)+1$

\bigskip

We will show that this case cannot occur. Consider the mapping from $\mathcal{B} \cup \mathcal{R}$ into $V$ defined as follows. To each $B_i$ we assign, using (\ref{four}), an element $u$ of $B_i$ which belongs to no other $B_{i'}$; if among the possible choices of $u$ for a given $B_i$ there is one which belongs to more than one of the sets $R_j$, we assign to $B_i$ such a $u$. Similarly, to each $R_j$ we assign an element $u$ of $R_j$ which belongs to no other $R_{j'}$, with priority to such $u$ which belongs to more than one of the sets $B_i$. As $|\mathcal{B} \cup \mathcal{R}| > |V|$, the mapping is not injective, so we can find some $u \in V$ which was assigned to some $B_i$ and to some $R_j$.

Assume w.l.o.g. that $b \ge 2(k-1)+1$. For the set $R_j$ just found, there is no vertex $v$ such that $|I_v| > 1$ and $J_v = \{j\}$; indeed, if there were such $v$ it would be given priority as the vertex assigned to $R_j$, over the actual assignment of $u$ which belongs to a unique $B_i$. It follows that we can partition $R_j$ into two sets:

\vspace{-5pt}

\begin{equation*}
S=\{v \in R_j: |I_v|=1 \}, \quad T=\{v \in R_j: |I_v|,|J_v| \ge 2 \}.
\end{equation*}

\noindent Write $|R_j| = k + \ell$, where $\ell \ge 0$, and $|S|=s$.

Due to the size of $R_j$, the difference between the two sides of (\ref{ten}) is at least $\ell$. Hence the first inequality in (\ref{thirteen}) holds with a slack of at least $\ell$.

For each $v \in T$, the difference between the two sides of (\ref{twelve}) is at least \[2(|I_v|+|J_v|-2) - (|I_v|+|J_v|-1) = |I_v|+|J_v|-3 \ge |I_v|-1.\] 
\noindent Since $|T|=k + \ell -s$, the second inequality in (\ref{thirteen}) holds with a slack of at least $\sum_{v \in T} |I_v| - (k + \ell - s)$. The sets $I_v$, $v \in R_j$, partition $[b]$, and therefore $\sum_{v \in T} |I_v| = b-s \ge 2(k-1)+1-s$, so the slack in the second inequality in (\ref{thirteen}) is at least $k - \ell -1$.

Adding up the two slacks, we obtain

\vspace{-15pt}

\begin{eqnarray*}
|V| & \ge & k(b+r) - br + k - 1\\
 & \ge & k(4(k-1)+1) - 2(k-1)(2(k-1)+1) + k - 1\\
 & = & 4(k-1) + 1,
\end{eqnarray*}

\noindent which contradicts our assumption on $|V|$.

\bigskip

It remains to address the possibility that the families $\mathcal{B}$, $\mathcal{R}$ associated with $G$ do not satisfy (\ref{four}) and (\ref{five}). In this case, by performing the steps indicated in the proof of Theorem~\ref{main}, we obtain modified families $\mathcal{B}'$, $\mathcal{R}'$ which do satisfy (\ref{four}) and (\ref{five}) as well as (\ref{one})--(\ref{three}). By the foregoing proof, $\mathcal{B}'$ and $\mathcal{R}'$ must be as described in Case~1 above, and the graph corresponding to them is isomorphic to some $G(X,Y,B,R)$. In particular, all sets in $\mathcal{B}' \cup \mathcal{R}'$ are of size $k$ exactly. It follows that in passing from $\mathcal{B}$, $\mathcal{R}$ to $\mathcal{B}'$, $\mathcal{R}'$, the step of adding a vertex to a set could never occur. Thus, the only steps performed were deletions of sets. Therefore the original graph $G$ contains a graph of the form $G(X,Y,B,R)$, and by edge-criticality they must coincide. $\square$

\section{Generalizations and reformulations}

\subsection{More than two colors}

It is natural to generalize the question treated here to $t$-colored graphs, i.e., simple graphs with edges colored in one of $t$ colors.

\begin{question}[Bucic et al.~\cite{BLLW}] \label{morecolors}
Given integers $k,t \ge 2$, what is the smallest possible number of vertices in a $t$-colored graph having the property that every vertex belongs to a monochromatic $k$-clique of each color?
\end{question}

Bucic et al. noted that their construction for $t=2$ on $4(k-1)$ vertices can be adapted to one for general $t$ using $2t(k-1)$ vertices. In fact, our more general construction in Section~3 can also be adapted as follows. Let $X_1,\ldots,X_t$ be $t$ disjoint sets of $2(k-1)$ vertices each. For any pair of colors $i,j$, take the complete bipartite graph on $(X_i,X_j)$ and color its edges $i$ or $j$ so that both color graphs are $(k-1)$-regular. In addition, for each color $i$, any two vertices in $X_i$ which have a common neighbor in the color $i$ graph (in any $X_j$, $j \ne i$) are joined by an edge colored $i$. This yields an edge-critical graph with the required property.

Regarding optimality, we note that the above construction is not optimal for $k=2$ and $t>2$. For example, $4$ vertices suffice for $k=2$, $t=3$. But it may be optimal for higher values of $k$. As observed by Bucic et al.~\cite{BLLW}, their proof of an asymptotic lower bound for the case of two colors extends to general $t$, yielding a lower bound of $(2t - o_k(1))k$ in Question~\ref{morecolors}. Unfortunately, it seems that our proof of the exact lower bound does not extend to general $t$.

While the above generalization looks interesting in its own right, the intended application of Bucic et al.~\cite{BLLW} suggests a different generalization. This will be explained in the following subsections.

\subsection{Partition of a box into sub-boxes}

A set of the form $A=A_1 \times \cdots \times A_d$, where $A_1,\ldots,A_d$ are finite sets with $|A_i| \ge 2$, is called a $d$-dimensional discrete box. A set of the form $B=B_1 \times \cdots \times B_d$, where $B_i \subseteq A_i$ for all $i \in [d]$, is a sub-box of $A$; it is said to be nontrivial if $\emptyset \ne B_i \ne A_i$ for all $i \in [d]$. It is easy to partition a $d$-dimensional discrete box into $2^d$ nontrivial sub-boxes, by cutting each $A_i$ into two parts. The following theorem answered a question of Kearnes and Kiss~\cite{KK}.

\begin{theorem}[Alon et al.~\cite{ABHK}] \label{book}
Let $A$ be a $d$-dimensional discrete box, and let $\{B^1,\ldots,B^m\}$ be a partition of $A$ into $m$ nontrivial sub-boxes. Then $m \ge 2^d$.
\end{theorem}

Instead of requiring the sub-boxes $B^1,\ldots,B^m$ to be nontrivial, one may equivalently require that every axis-parallel line (i.e., set of the form $\{(x_1,\ldots,x_d) \in A: x_j = a_j \,\,\forall j \in [d] \setminus \{i\}\}$) intersects at least two of them. This observation led Bucic et al.~\cite{BLLW} to consider families of sub-boxes $\{B^1,\ldots,B^m\}$ with the $k$-piercing property, namely: every axis-parallel line intersects at least $k$ sub-boxes in the family. Generalizing the question of Kearnes and Kiss, they asked the following.

\begin{question}[Bucic et al.~\cite{BLLW}] \label{piercing}
Let $d \ge 1$ and $k \ge 2$ be integers, and let $A=A_1 \times \cdots \times A_d$ be a $d$-dimensional discrete box with all $|A_i|$ sufficiently large. What is the smallest possible number $m$ of sub-boxes in a partition $\{B^1,\ldots,B^m\}$ of $A$ having the $k$-piercing property?
\end{question}

They denoted the answer to Question~\ref{piercing} by $p_{\mathrm{box}}(d,k)$. The case $k=2$ is solved by Theorem~\ref{book}: $p_{\mathrm{box}}(d,2)=2^d$. For larger $k$, it is natural to consider first the $2$-dimensional case ($d=1$ is trivial). Here, cutting each $A_i$ into $k$ parts gives a construction with $m=k^2$ sub-boxes. But Bucic et al.~\cite{BLLW} showed that in fact $m=4(k-1)$ is enough. Their construction is illustrated in Figure~1.

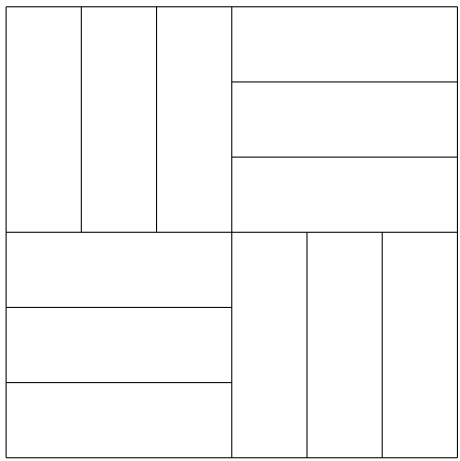
\begin{figure}[htbp]
\begin{center}
\begin{tikzpicture}
\draw (-3,3) -- (3,3);
\draw (0,2) -- (3,2);
\draw (0,1) -- (3,1);
\draw (-3,0) -- (3,0);
\draw (-3,-1) -- (0,-1);
\draw (-3,-2) -- (0,-2);
\draw (-3,-3) -- (3,-3);
\draw (-3,-3) -- (-3,3);
\draw (-2,0) -- (-2,3);
\draw (-1,0) -- (-1,3);
\draw (0,-3) -- (0,3);
\draw (1,-3) -- (1,0);
\draw (2,-3) -- (2,0);
\draw (3,-3) -- (3,3);
\end{tikzpicture}

\caption{A $k$-piercing partition of a $2$-dimensional box, showing that $p_{\mathrm{box}}(2,k) \le 4(k-1)$. Each quarter of the box consists of $k-1$ parallel sub-boxes.}

\end{center}
\end{figure}

Bucic et al. conjectured that this construction is optimal, that is, $p_{\mathrm{box}}(2,k)=4(k-1)$. They observed that this is the case if one restricts attention to sub-boxes which are bricks, i.e., products of intervals. In an attempt to prove optimality among partitions into general sub-boxes, they associated with any such partition of a $2$-dimensional box a $2$-colored graph as follows: the vertices are the sub-boxes in the partition, and two sub-boxes are joined by a blue (resp. red) edge if there is a horizontal (resp. vertical) line which intersects both of them. The $k$-piercing property implies that every vertex belongs to a monochromatic $k$-clique of each color. This led them to ask for the minimum number of vertices in such a graph. Note that the $2$-colored graph with $4(k-1)$ vertices constructed by them (and presented in the introduction) corresponds to the partition shown in Figure~1.

The asymptotic lower bound that Bucic et al.~\cite{BLLW} obtained for the question about $2$-colored graphs enabled them to deduce that $p_{\mathrm{box}}(2,k) \ge (4 - o_k(1))k$. Our full solution of the question (Theorem~\ref{main}) allows us to confirm their conjecture: $p_{\mathrm{box}}(2,k)=4(k-1)$. In fact, since the reduction described in the previous paragraph does not depend on the sub-boxes being a covering of the given box, but only on their disjointness, we have the following more general statement.

\begin{corollary} \label{partition}
Let $k \ge 2$ be an integer, let $A$ be a $2$-dimensional discrete box, and let $\{B^1,\ldots,B^m\}$ be a family of $m$ disjoint sub-boxes of $A$ having the $k$-piercing property. Then $m \ge 4(k-1)$.
\end{corollary}

The question of determining $p_{\mathrm{box}}(d,k)$ when both $d$ and $k$ are greater than $2$ remains wide open. Bucic et al.~\cite{BLLW} attempted a reduction to colored graphs similar to the above, but it led to a less natural and less tractable question than in the case $d=2$. Their best bounds for general $d$ and $k$ are of the form $e^{\Omega(\sqrt{d})}k \le p_{\mathrm{box}}(d,k) \le 15^{d/2}k$ (of course, when $k$ is small relative to $d$, the bounds $2^d = p_{\mathrm{box}}(d,2) \le p_{\mathrm{box}}(d,k) \le k^d$ may be better).

\subsection{Decomposition of a bipartite graph into complete bipartite subgraphs}

A well-studied parameter of a graph $G=(V,E)$ is the minimum number of edge-disjoint complete bipartite subgraphs of $G$ which cover the edge-set $E$. The best known result is that of Graham and Pollak~\cite{GP}, saying that any such decomposition of the complete graph $G=K_n$ must consist of at least $n-1$ complete bipartite subgraphs. For more general results about decomposition of an arbitrary graph $G$, see e.g. Kratzke et al.~\cite{KRW}. The case when $G$ itself is complete bipartite is of course uninteresting, because there is a decomposition into one subgraph. But it becomes interesting under some constraints on the decomposition, as we will see below.

A $2$-dimensional discrete box $A=A_1 \times A_2$ (discussed in the previous subsection) may be viewed as the edge-set of a complete bipartite graph on $(A_1,A_2)$. A partition of $A$ into sub-boxes is then a decomposition of a complete bipartite graph into complete bipartite subgraphs. We can restate Corollary~\ref{partition} from this point of view, as follows.

\begin{corollary} \label{decomposition2k}
Let $k \ge 2$ be an integer. Let $G=(A_1,A_2,E)$ be a bipartite graph, and let $\{ G^i = (B_1^i,B_2^i,E^i) \}_{i \in [m]}$ be a decomposition of $G$ into $m$ complete bipartite subgraphs. Assume that every vertex in $A_1$ (resp. $A_2$) belongs to $B_1^i$ (resp. $B_2^i$) for at least $k$ values of $i \in [m]$. Then $m \ge 4(k-1)$.
\end{corollary}

Proposition~\ref{prop} may also be reformulated in this terminology, as follows.

\begin{corollary} \label{decomposition2}
Let $G=(A_1,A_2,E)$ be a complete bipartite graph, and let $\{ G^i = (B_1^i,B_2^i,E^i) \}_{i \in [m]}$ be a decomposition of $G$ into $m$ complete bipartite subgraphs. Assume that for every vertex $x$ in $A_1$ (resp. $A_2$) there is $i \in [m]$ such that $B_1^i = \{x\}$ (resp. $B_2^i = \{x\}$). Then $m \ge |A_1|+|A_2|-1$.
\end{corollary}

Indeed, Tverberg's~\cite{T} proof of Graham and Pollak's theorem inspired the proof of Proposition~\ref{prop}.

This point of view on partition problems for $2$-dimensional discrete boxes suggests a generalization to higher dimensions expressed in terms of $d$-partite hypergraphs. In particular, the following question asks for a $d$-partite version of Corollary~\ref{decomposition2k}.

\begin{question} \label{decompositiondk}
Let $d,k \ge 2$ be integers. Let $H=(A_1,\ldots,A_d,E)$ be a complete $d$-partite hypergraph, and let $\{ H^i = (B_1^i,\ldots,B_d^i,E^i) \}_{i \in [m]}$ be a decomposition of $H$ into $m$ complete $d$-partite subhypergraphs. Assume that for every $\ell \in [d]$ and for every $(d-1)$-tuple of vertices $x_j \in A_j$, $j \in [d] \setminus \{\ell\}$, there are at least $k$ values of $i \in [m]$ such that $x_j \in B_j^i$ for all $j \in [d] \setminus \{\ell\}$. If all $|A_j|$ are sufficiently large, what is the smallest possible number $m$ of subhypergraphs in such a decomposition?
\end{question}

This is a reformulation of Question~\ref{piercing}, so the answer is the same $p_{\mathrm{box}}(d,k)$ investigated by Bucic et al.~\cite{BLLW}. Hopefully, this interpretation of the question may suggest a useful approach, but we were unable to extend the methods of this paper to handle it.

\end{document}